\documentstyle[11pt,amssymb,fullpage]{article}
\newtheorem{theo}{Theorem}[section]
\newtheorem{coro}[theo]{Corollary}

\newtheorem{lem}[theo]{Lemma}

\def\SBIMSMark#1#2#3{
 \font\SBF=cmss10 at 10 true pt
 \font\SBI=cmssi10 at 10 true pt
 \setbox0=\hbox{\SBF Stony Brook IMS Preprint \##1}
 \setbox2=\hbox to \wd0{\hfil \SBI #2}
 \setbox4=\hbox to \wd0{\hfil \SBI #3}
 \setbox6=\hbox to \wd0{\hss
             \vbox{\hsize=\wd0 \parskip=0pt \baselineskip=10 true pt
                   \copy0 \break%
                   \copy2 \break% 
                   \copy4 \break}}
 \dimen0=\ht6   \advance\dimen0 by \vsize \advance\dimen0 by 8 true pt
                \advance\dimen0 by -\pagetotal
 \dimen2=\hsize \advance\dimen2 by .25 true in
  \newread\jref
  \openin\jref=publishd.tex
  \ifeof\jref\setbox0=\hbox to 0pt{}
  \else 
     \setbox0=\hbox to 3.1 true in{
                \vbox to \ht6{\hsize=3 true in \parskip=0pt  \noindent  
                \input publishd.tex 
                \vfill}}
  \fi
  \closein\jref
  \ht0=0pt \dp0=0pt
 \ht6=0pt \dp6=0pt
 \setbox8=\vbox to \dimen0{\vfill \hbox to \dimen2{\copy0 \hss \copy6}}
 \ht8=0pt \dp8=0pt \wd8=0pt
 \copy8
 \message{*** Stony Brook IMS Preprint #1, #2 ***}
}

\begin{document}

\begin{center}
{\Large \bf Non-accessible Critical Points of Cremer Polynomials}\\
\vspace{2 cm}
Jan Kiwi\\
Mathematics  Deparment\\
S.U.N.Y.  Stony Brook\\
\end{center}

\begin{abstract}
 It is shown  that a polynomial with a Cremer periodic point has a
non-accessible critical point in its Julia set provided that the
Cremer periodic point is approximated by small cycles. 
\end{abstract}

\SBIMSMark{1995/2}{February 1995}{}
\thispagestyle{empty}

\section{Introduction.}

Cremer fixed points occur as fixed points of polynomial
maps from the complex plane  onto itself.
Their presence forces the Julia set not to be locally connected 
\cite{su:geo} \cite{lyu:lec}.
When local connectivity fails it is of interest to 
determine which points in the Julia set are accessible
from the complement of the  Julia set and which are not (if any).
\medskip

Repelling periodic points are always accessible \cite{hu:yoc} \cite{el:acc}
 and 
it `often' happens that Cremer fixed points are
`surrounded' by repelling cycles. Let us be more 
precise, if every neighborhood of a Cremer fixed point
contains infinitely many repelling cycles, we
say that the Cremer fixed point is approximated by small
cycles. For quadratic polynomials,
J.-C. Yoccoz has shown that this is always the case \cite{yo:linea}. For 
polynomials of higher degree it is an open problem.
 Some results and a more detailed discussion about
this problem can be found in  \cite{milnor:section8}  and \cite{pm:cremer}.
\medskip

In this paper we show how the accessibility of the small cycles
implies the non-accessibility of a critical point.
 
\begin{theo}
Let $P$ be a polynomial with a Cremer fixed point $\hat{z}$
that is approximated by small cycles. Then there exists
a critical point which is not accessible 
from the complement of the Julia set.
\end{theo}

We conclude that the critical
point of a quadratic polynomial with a Cremer periodic point  is not
accessible, thus answering a question posed by J. Milnor in
\cite{bi:problem}.

\begin{coro}
The critical point of a quadratic polynomial with a Cremer
periodic point is not accessible from the complement of the 
Julia set.
\end{coro}

Under the assumption that the Julia set 
is connected, each point that is accessible from
the complement of the filled Julia set is
the landing point of some external ray.
In this case, we prove a stronger statement by showing the 
existence of a non-accessible critical value with the following property.
An external ray (if any) that accumulates at
this critical value also accumulates at the 
Cremer fixed point.

\begin{theo}
Let $P$ be a polynomial with connected  Julia set $J_P$
and with a Cremer fixed point $\hat{z}$ that is approximated
by small cycles. Then there exists a 
 critical value $v$ which is not accessible from  ${\Bbb C} \smallsetminus J_P$ such
that: if $v$ belongs to the closure $\overline{R}$ of some
external ray $R$, 
then $\hat{z}$ also belongs to $\overline{R}$.
\end{theo}

For some quadratic polynomials our result can be de\-duced from a
stronger one of A. Douady (compare \cite{so:mas}). More precisely, for
a generic quadratic polynomial with a Cremer fixed point,  Douady
showed that there is an external ray that does not land; furthermore,
this external ray accumulates in a set containing the Cremer fixed
point and its preimage.  It follows that for such quadratic
polynomials every external ray which accumulates at the critical value
also accumulates at the Cremer fixed point (see Corollary 3.5, Section 3).
\medskip

Recently,  R. P\'erez-Marco \cite{pm:hedge} has  shown that rational
functions with Cremer periodic points always have some non-accessible points in
their Julia set. For a `large' family of infinitely renormalizable
quadratic polynomials, he proves the non-accessibility of the critical
point; and for a generic quadratic polynomial with a Cremer fixed point,
P\'erez-Marco shows that the critical point is not accessible.
\bigskip

%\newpage

\noindent
This paper is organized as follows:
\smallskip

In Section 2, we briefly summarize some definitions and  results
from polynomial dynamics. For general background in complex
dynamics see \cite{milnor:lec} and \cite{lyu:lec}.
\smallskip

In Section 3, we prove some preliminary results. Also, we show how one
of them leads to a new proof (for polynomials with connected Julia
set) of a well known result by A. Douady \cite{bl:lec} and M. Shishikura
\cite{sh:non}. That is, {\it for a polynomial of degree $d$,
 the total number of Cremer cycles and cycles of  bounded Fatou components
is strictly less than $d$}
 (Corollary 3.4).
\smallskip

In Section 4, we prove  our main results.
 The proof of Theorem 1.3 is by contradiction and the idea is
a refinement of the following line of argument. Under the assumption
that certain critical values  are accessible, we construct
an open region around the Cremer fixed point containing at most
finitely many cycles. This region is obtained by cutting the complex
plane along an appropriate collection of external rays. The
correct choice of these external rays is based
on the results from Section 3. In order to prove Theorem 1.1 
in the disconnected case we use some
results about polynomial like mappings and then we  apply Theorem~1.3
to find a non-accessible critical point. 
\bigskip

\noindent
{\bf Acknowledgments:}
I am grateful  to  John Milnor for 
many interesting and inspiring discussions. 
I would like to thank Misha Lyubich and Alfredo Poirier
for many helpful conversations.
The author is partially supported by a {\it ``Beca Presidente de la Republica'' (Chile)}. 

\section{Cremer Points and Accessible Points.}

\indent
Consider a polynomial $P: {\Bbb C} \rightarrow {\Bbb C}$ of degree $d
\geq 2$ with a fixed point $\hat{z}$ of {\it multiplier} $\lambda =
P^{\prime} (\hat{z})$. The fixed point is called {\it repelling} if $
|\lambda| > 1 $, {\it attracting} if $| \lambda | < 1 $ and {\it
parabolic} if $\lambda$ is a root of unity.  
\smallskip

For $\lambda = e^{ 2 \pi
i \theta}$ where $\theta$ is irrational, we have to distinguish
between two possibilities. If the fixed point $\hat{z}$ belongs to the
Julia set $J_P$ we say that $\hat{z}$ is a {\it Cremer fixed
point}. Otherwise, if $\hat{z}$ lies in the Fatou set $F_P$, call
$\hat{z}$ a {\it Siegel fixed point}.  
\medskip

For a generic choice of $\theta$, H. Cremer \cite{cr:zentru} showed
that there exists a sequence of periodic points converging to
$\hat{z}$, in this case $\hat{z}$ belongs to the Julia set.  In
contrast, for a full Lebesgue measure set of angles $\theta$,
C.L.~Siegel \cite{si:linea} proved that $\hat{z}$ lies in the Fatou
set.  We refer the reader to \cite{milnor:section8} and \cite{pm:cremer}
for a more detailed exposition about this dichotomy.
\medskip

A modification of Cremer's proof (compare \cite{milnor:section8})
shows that for a generic choice of $\theta$, every neighborhood of a
Cremer fixed point $\hat{z}$ contains infinitely many cycles.  We say
that $\hat{z}$ {\it is approximated by small cycles}. According to
Yoccoz \cite{yo:linea} quadratic Cremer fixed points are always
approximated by small cycles. It is not known whether this is true for
polynomials of arbitrary degree. For some interesting results about
this problem  see \cite{pm:cremer}.
\bigskip

A point $z$ in the Julia set $J_P$ is called {\it accessible} 
from a Fatou component $V$ if there
exists a path $\gamma$ contained in $V$ that ends at $z$, i.e. 
$$\gamma : [0,1) \rightarrow V$$ and
$\gamma ( r ) \rightarrow z$ as $r \rightarrow 1$.  Accessibility is a
property completely invariant under the dynamics, $z$ is accessible 
from $V$ if
and only if $P(z)$ is accessible from $P(V)$. When $z$ is not accessible
from any Fatou component sometimes we just say that {\it $z$ is a non-accessible point}.
\medskip

Let us now assume that $K_P$ is connected. For simplicity
we also assume that $P$ is a monic polynomial. The complement
of the closed unit disc $\overline{\Delta}$ is conformally isomorphic to the 
complement of  $K_P$. Moreover, there exists a unique conformal
isomorphism: 
 $$\phi :{\Bbb C} \smallsetminus \overline{\Delta} \rightarrow 
{\Bbb C} \smallsetminus K_P$$
such that  $$P \circ \phi ( z ) = \phi ( z^d ),$$
and that is asymptotic to identity near infinity (see \cite{milnor:lec}).
\medskip

The question of which points in $J_P$ are accessible
from ${\Bbb C} \smallsetminus K_P$ is closely
related to the boundary behavior of the map $\phi$.  We say that 
$$R_t = \phi ( \{ r e^{2 \pi i t} : r \in (1, \infty) \} )$$ is the
 {\it external ray} with angle $t \in {\Bbb R}/ {\Bbb Z}$.  An
 external ray $R_t$ {\it lands at} $z$ if $\phi ( r e^{2 \pi i t} )
 \rightarrow z$ as $r \rightarrow 1$.  Therefore, the landing points
 of external rays are accessible from ${\Bbb C} \smallsetminus K_P$.
 A result of Lindel\"of assures that the converse is also true, every
 point accessible from ${\Bbb C} \smallsetminus K_P$ is the landing
 point of some external ray (compare \cite{ru:real}).
\medskip
 
The external ray $R_t$ is called a {\it periodic ray} if $P^{\circ n} ( R_t ) =
R_t$ for some positive integer $n$.  In particular, if $P(R_t)= R_t$
we say that the external ray $R_t$ is a {\it fixed ray}.  Periodic
rays always land, moreover their landing points are parabolic or
repelling periodic points.  Conversely, parabolic and repelling
periodic points are the landing points of at least one periodic ray.
(see \cite{hu:yoc} \cite{milnor:lec})
\medskip

If an external ray $R$ lands, then the closure $\overline{R}$ is the
union of the external ray and its landing point. Otherwise, if $R$
fails to land,  $\overline{R}$ is the union of the external ray and
a non-trivial connected subset of the Julia set. 
For short, the closure  $\overline{R}$  of an external ray $R$ is called
 a {\it closed ray}.

\section{Separation Lemmas.}

For polynomials with connected Julia sets, A. Poirier pointed out that
the union  of all closed fixed rays `separates' invariant Fatou
components and Cremer fixed points. This is stated  as a
corollary of a result by L. Goldberg and J. Milnor (compare with
Theorem 3.3 and Corollary 3.5 in
\cite{gol-mil:por}). We go one step further and show that after
cutting the complex plane along an appropriate collection of external
rays, each Cremer periodic point and periodic Fatou component is
`separated' from the others:

\begin{lem}
Let $P$ be a polynomial with connected Julia set. Then there
exists a finite union ${\cal C}$ of closed periodic rays such that:
\smallskip

  $P ({\cal C}) = {\cal C} $.

  Each component of ${\Bbb C} \smallsetminus {\cal C}$ contains at most
one Cremer periodic point or periodic Fatou component.
\end{lem}

Recall that periodic rays always land, hence the corresponding closed
ray is the union of the periodic ray and its landing point. In this
section, the closed rays taken into consideration are always
the closure of external rays that land.
\medskip

Lemma  3.1 is an immediate consequence of the following (compare
with Theorem 3.3 in \cite{gol-mil:por}):

\begin{theo}[Goldberg and Milnor]
Let $P$ be a polynomial with connected Julia set. Denote by
${\cal C}$ the union  of all closed fixed rays:
$${\cal C} =  \overline{R}_{0} \cup  \overline{R}_{1/d-1} \cup \cdots
\cup  \overline{R}_{d-2/d-1}$$
If $U$ is a connected component of   ${\Bbb C} \smallsetminus {\cal C}$, 
then one
and only one of the following holds:
\smallskip

 $U$ contains exactly one invariant Fatou component.

 $U$ contains exactly one fixed point which belongs to
the Julia set. 
\end{theo}

\newpage

\noindent
{\bf Proof of Lemma 3.1:} If $P$ has no Cremer periodic points and no
bounded periodic Fatou components there is nothing to prove.
Otherwise, let $m$ be the minimum common multiple of the following
list of integers: periods of periodic Fatou components, and periods of
Cremer periodic points.  The number of non-repelling cycles is finite,
hence this list is also finite and $m$ is a well defined integer.
\smallskip

Periodic Fatou components and Cremer periodic points of $P$ and
$P^{\circ m}$ are the same, only their periods differ. In fact, we
have chosen $m$ in such a way that the periodic Fatou components are
invariant under $P^{\circ m}$ and that the Cremer periodic points are fixed
points in its Julia set.
\smallskip

Now consider the union ${\cal C}$ of closed rays
fixed by $P^{\circ m}$:
$${\cal C}=  \overline{R_0} \cup \cdots \cup \overline{R}_{d^m -2/d^m -1}.$$
Since fixed rays land at parabolic or repelling points, this
set ${\cal C}$ is disjoint from bounded Fatou components and 
Cremer points.
Applying the previous Theorem to $P^{\circ m}$, it follows that
the set ${\cal C}$ has the desired properties. \hfill $\Box$ 
\bigskip

Let us introduce some notation. If ${\cal C}$ is the union of
finitely many
closed external rays and $z$ is a point that lies in the complement of
${\cal C}$, denote by $U_ {\cal C} ( z )$ the connected component of
${\Bbb C} \smallsetminus {\cal C}$ that contains $z$. 
In a similar way, if $V$ is a
connected set contained in the complement of ${\cal C}$, denote by $U_
{\cal C} ( V )$ the connected component of ${\Bbb C} \smallsetminus {\cal C}$
 that
contains $V$.
\medskip

A refinement of the previous Lemma shows that  cutting along
a larger collection of external rays, critical points and Cremer
points that are left `together',  stay `together' under 
forward iterations. The corresponding  statement  for periodic
Fatou components also holds.

\begin{lem}
Let $P$ be a polynomial with connected Julia set. Then 
there exists a finite union ${\cal C}$ of closed rays
such that:
\smallskip

{\rm (1)}  $P ({\cal C}) \subset {\cal C} $.

{\rm (2)}  Each component of ${\Bbb C}-{\cal C}$ contains at most
one Cremer periodic point or periodic Fatou component.

{\rm (3)} If $\hat{z}$ is a Cremer periodic point and $c$ is a
critical point such that $U_{\cal C} (c) = U_ {\cal C} ( \hat{z} )$,
then $U_{\cal C} ( P^{\circ n} ( c )) = U_ {\cal C}
(P^{\circ n} ( \hat{z} ) )$ for all $n \geq 1$. Moreover, $c \in J_P$ and
$c$ does not belong to the closure of any bounded Fatou component.

{\rm (4)} If $V$ is a bounded periodic Fatou component and 
$c$ is a critical point such that $U_{\cal C} (c) = U_ {\cal C} ( V )$,
then $U_{\cal C} (P^{\circ n} ( c )) =  U_ {\cal C} (P^{\circ n} ( V ) )$
 for all $n \geq 1$. 
\end{lem}

\noindent
{\bf Remark:} The proof of Corollary 3.4 will show that 
statements (3) and (4) of the previous Lemma are not vacuous
in the following  sense. Given  a  Cremer periodic  orbit  
$z_0 = P ( z_{n-1} ) , \dots, z_{n-1} = P (  z_{n-2}  )$,  then
for some $z_i$  there exists a critical point $c$ such that 
 $U_{\cal C} (c) = U_ {\cal C} ( z_i )$.  The corresponding assertion
is also true for bounded periodic Fatou  components. \medskip

\noindent
{\bf Proof:} Let ${\cal C}_0$ be the union of closed  rays
from Lemma 3.1. For $k \geq 1$ define inductively
$ {\cal C}_k =  P^{-1} ({\cal C}_{k-1}).$
\smallskip

We are going to show that conditions (1) through (4) hold for 
${\cal C}_k$ when $k$ is large enough. Before let us list some 
properties of the sets  ${\cal C}_k$ and their complementary
regions:
\smallskip

%\newpage

(a) ${\cal C}_k \subset   {\cal C}_{k+1}$ for $k \geq 0$,

(b) $P ( {\cal C}_k ) \subset  {\cal C}_k$ for $k \geq 0$,

(c) $U_{{\cal C}_{k+1}}(z) \subset U_{{\cal C}_k}(z)$ for $k \geq 0$ and 
$z \in {\Bbb C}-{\cal C}_{k+1}$,

(d) $P^{\circ n}(U_{{\cal C}_{k+n}} (z)) = 
U_{{\cal C}_k} (P^{\circ n} (z))$ for $k \geq 0$, $n \geq 1$ and 
$z \in {\Bbb C}-{\cal C}_{k+n}$.
\smallskip

All of the above can be deduced from the forward invariance of $ {\cal C}_0$
 ($ P({\cal C}_0) = {\cal C}_0$).
\medskip

From Lemma 3.1 and properties (b) and (c) it follows that conditions
(1) and (2) hold for all ${\cal C}_k$.
\medskip

{\it Claim 1:} Given a critical point $c$ and a Cremer periodic
point $\hat{z}$ there exists $k(c, \hat{z}) \geq 0$ such that:

\noindent
($\ast$) if $U_{{\cal C}_l}(c) = U_{{\cal C}_l} (\hat{z})$, then 
 $U_{{\cal C}_l} ( P^{\circ n} ( c )) = U_{{\cal C}_l}
(P^{\circ n} ( \hat{z} ) )$ for all $n \geq 1$; 

\noindent
holds for all $l \geq k(c, \hat{z})$.
\medskip

{\it Proof of  Claim 1:} 
In the case that ($\ast$) holds for all $l \geq 0$, 
of course we take  $k(c, \hat{z}) = 0$.
Otherwise, for some  $ l_0 \geq 0$ we have that  $U_{{\cal C}_{l_0}}(c) =
U_{{\cal C}_{l_0}} (\hat{z})$
and 
$U_{{\cal C}_{l_0}} ( P^{\circ n} ( c )) \neq U_{{\cal C}_{l_0}} (P^{\circ n}
( \hat{z} ) )$ for some $n \geq 1$.
In this case let  $k(c, \hat{z}) = l_0 + n$, it
follows from (c) and (d) that 
 $U_{{\cal C}_l}(c) \neq U_{{\cal C}_l} (\hat{z})$ for all 
$l \geq  k(c, \hat{z})$.\medskip

We omit the identical proof of the corresponding claim for periodic
Fatou components:\medskip

{\it Claim 2:} Given a critical point $c$ and a bounded periodic
Fatou component $V$ there exists $k(c, V) \geq 0$ such that: 

\noindent
($\ast$) if $U_{{\cal C}_l}(c) = U_{{\cal C}_l} (V)$, then 
 $U_{{\cal C}_l} ( P^{\circ n} ( c )) = U_{{\cal C}_l}
(P^{\circ n} ( V ) )$ for all $n \geq 1$;

\noindent
holds for all  $l \geq k(c, V)$.\medskip

Let $N$ be the maximum of the numbers $k(c, \hat{z})$ and $k(c, V)$,
 where all possible pairs $(c, \hat{z})$ and $(c, V)$ are taken into
 consideration. 
Since the number of critical points and non-repelling cycles is
finite, the total  number of all such pairs is also finite and
$N$ is a well defined integer.
\medskip

As mentioned before Claim 1,  ${\cal C}_N$ satisfies conditions
(1) and (2). From Claim 2, it follows that (4) holds for  ${\cal C}_N$.
It remains to complete the proof that (3) also holds for  ${\cal C}_N$.
\smallskip

 If $c$ is a critical point and $\hat{z}$ is a
 Cremer point such that $U_{{\cal C}_N}(c) = U_{{\cal C}_N} (\hat{z})$,
 then 
$$U_{{\cal C}_N} ( P^{\circ n} ( c )) = U_{{\cal C}_N} (P^{\circ
 n} ( \hat{z} ) ) \neq U_{{\cal C}_N} (V)$$ 
for all $n \geq 1$ and for all bounded periodic Fatou
 components $V$. According to D. Sullivan \cite{su:nowand},
 every point in the closure of a bounded Fatou component
 is mapped after a finite number of iterates to the
closure of a bounded periodic Fatou
 component. It follows that $c$ belongs to $J_P$ and $c$ does not 
belong to the closure of any bounded Fatou component.
\hfill $\Box$

\bigskip
According to Douady \cite{bl:lec} and Shishikura \cite{sh:non},
 the  total number
of Cremer cycles and cycles of bounded Fatou components of $P$ is strictly less
than $d = \deg P$. For polynomials with connected Julia sets, 
Lemma 3.1 leads to a new proof of this result.

\begin{coro}
Let $P$ be a polynomial with connected Julia set.
Then the total number of Cremer cycles and cycles of bounded Fatou components
is  less than or equal to the number of distinct  critical points.
\end{coro}

\noindent
{\bf Proof:} The idea is to show that to each Cremer cycle 
corresponds at least one critical point and to each cycle 
of bounded Fatou components also corresponds at least one critical
point.
\smallskip

Apply Lemma 3.1 to obtain a finite union ${\cal C}$ of closed 
external rays. Consider a Cremer cycle
 $z_0 = P(z_{n-1}),
\dots, z_{n-1} = P(z_{n-2})$.
 We claim that there exists at 
least one critical point in 
$U_{\cal C} (z_0) \cup \cdots
\cup U_{\cal C} (z_{n-1})$.
 In order to prove this claim let ${\cal C}_1 = P^{-1} ({\cal C})$
and for $k \geq 2$ 
define inductively ${\cal C}_k = P^{-1} ( {\cal C}_{k-1})$.
Suppose that $U_{\cal C} (z_0) \cup \cdots
\cup U_{\cal C} (z_{n-1})$ contains no critical points.
In particular, $U_{{\cal C}_n} (z_0) \cup \cdots
\cup U_{{\cal C}_1} (z_{n-1})$ contains no critical point.
It follows that the proper holomorphic map
$$ P^{\circ n} |_{U_{{\cal C}_n} (z_0)}: {U_{{\cal C}_n} (z_0)}
\rightarrow U_{\cal C} (z_0)$$
has no critical points. Therefore,  
$ P^{\circ n} |_{U_{{\cal C}_n} (z_0)}$ is a conformal isomorphism
of simply connected domains. By Schwartz Lemma,
$ {U_{{\cal C}_n} (z_0)} = U_{\cal C} (z_0) $ and 
$ P^{\circ n} |_{U_{{\cal C}_n} (z_0)}$ is conformally conjugate
to a rotation of the unit disc, which contradicts the fact
that $K_P$ is bounded. Hence, 
$U_{\cal C} (z_0) \cup \cdots
\cup U_{\cal C} (z_{n-1})$ contains at least one critical point. 
\medskip

For a cycle   $V_0 = P(V_{n-1}), \dots, V_{n-1} = P(V_{n-2})$
of attracting or parabolic basins 
 it is well known that
$U_{\cal C} (V_0) \cup \cdots
\cup U_{\cal C} (V_{n-1})$ contains at least one critical point
(see \cite{milnor:lec}).
For a cycle of Siegel discs, an identical argument then the 
one used for Cremer cycles shows that $U_{\cal C} (V_0) \cup \cdots
\cup U_{\cal C} (V_{n-1})$ contains at least one critical point.
The proof follows from the classification of periodic Fatou components.
\hfill $\Box$
\bigskip

In order to prove the remark after Lemma 3.3, repeat the proof
of the previous Corollary but this time start with the
finite union ${\cal C}$ of closed rays obtained from Lemma 3.3.\medskip

The next result shows how the existence of an external ray that
accumulates at a Cremer fixed point and one of its preimages implies
the existence of a non-accessible critical value. Moreover, the
non-accessible critical value has the property that if an external ray
accumulates at it, this external ray also accumulates at the Cremer
fixed point.

\begin{coro}
Let $P$ be a polynomial with connected $K_P$ and with a Cremer
fixed point $\hat{z}$.  If there exists an external ray
$R_0$ such that $\{ \hat{z} , z \} \subset \overline{R_0}$ for
some $z \in P^{-1} (\hat{z}) \smallsetminus \{ \hat{z} \}$,
then there exists a non-accessible  critical value $v$ with the
following property: every external ray $R$ that accumulates
at $v$ also accumulates at $\hat{z}$.
\end{coro}

\noindent
{\bf Outline of the Proof:} Apply Lemma 3.3 to obtain a finite union ${\cal C}$
of closed rays. 
Let ${\cal C}_1 = P^{-1}({\cal C})$ and 
notice that $R_0 \subset U_{{\cal C}_1} (\hat{z})$.
Call $c_1 , \dots , c_k$ the critical points that are contained
in $ U_{{\cal C}_1} (\hat{z})$. It follows that 
their critical values $v_1 , \dots, v_k$  
are contained in  $ U_{{\cal C}} (\hat{z})$, and  the map
$$ P:  U_{{\cal C}_1} (\hat{z}) 
\smallsetminus P^{-1} (\{ v_1, \dots , v_k \}) \rightarrow  
U_{{\cal C}} (\hat{z}) \smallsetminus \{ v_1, \dots , v_k \}$$
is a regular covering.
\medskip

Take a small disc $D$ centered at $\hat{z}$. Call $D_0$  the preimage of $D$
that contains  $\hat{z}$ and denote by $D_1 , \dots , D_l$ 
the other preimages of $D$ contained in
$ U_{{\cal C}_1} (\hat{z})$. It follows that there
exists a closed segment  $\alpha$ of $R_0$ 
contained in  $ U_{{\cal C}_1} (\hat{z}) \smallsetminus D_0 \cup \cdots \cup
D_l$
such that $\alpha$ joins the boundary $\partial D_0$ of $D_0$
 with the boundary $\partial D_i$ of $D_i$ for some $i \neq 0$.
Also, we can assume that $\alpha$ intersects $\partial D_0 \cup \partial D_i$
 in 
exactly 2 points, which are the endpoints of $\alpha$. 
\medskip

Observe that $P(\alpha)$ is contained in  $U_{{\cal C}} (\hat{z}) 
\smallsetminus D$ and intersects $\partial D$ in exactly 2 points.
The set $\partial D  \cup P(\alpha)$ separates the plane into 3 components.
That is, an unbounded component, the disc $D$, and a  bounded component
$\tilde{D}$. Call $\beta$ the boundary of the Jordan domain 
$\tilde{D} \subset  U_{\cal C} (\hat{z})$. The curve $\beta$ is 
the union of a segment of $\partial D$ with $P(\alpha)$. Since a lift
of $\beta$ is not a closed curve, $\tilde{D}$ contains a critical value
$v$.
By Lemma 3.3 (3) the 
critical value $v$  does not belong to the closure of any bounded
Fatou component. 
 If an external ray $R$ accumulates at $v$, then $R$ intersects
$\partial D$. Taking $D$ arbitrarily small, it follows that
$v$ has the desired properties.
\hfill $\Box$

\section{Non-accessible Points.}

In this section first we prove the following:

\medskip
\noindent
{\bf Theorem 1.3} {\it Let $P$ be a polynomial
with connected Julia set $J_P$ and with a Cremer
fixed point $\hat{z}$ that is
approximated by small cycles. Then there exists a 
critical value $v$  which is not accessible
from ${\Bbb C} \smallsetminus J_P$ such that: if $v \in \overline{R}$ for
some external ray $R$, then $\hat{z} \in \overline{R}$ }.
\medskip

\noindent
{\bf Proof of Theorem 1.3}: Start applying Lemma 3.3 
to obtain a finite union ${\cal C}$ of closed rays.
Notice that all the critical points $c_1 , \dots, c_k$  contained
in $U_{{\cal C}}(\hat{z})$
belong to the Julia set $J_P$. Thus, their 
critical values $ v_1 = P (c_1), \dots, v_k = P (c_k)$  are also in $ J_P$.
These critical values are not in the closure of any bounded Fatou
component. Therefore, it is enough to show that one of them is
not accessible from ${\Bbb C} \smallsetminus K_P$ and that
every external ray which accumulates at this critical value 
also accumulates at the Cremer fixed point.
\smallskip

Let us first consider the case in which the Cremer fixed point
is not a critical value.
The proof is by contradiction. Assume that there exists
external rays $R_{t_1}, \dots, R_{t_k}$ such that 
$v_i \in \overline{R_{t_i}}$ and $\hat{z} \notin  \overline{R_{t_i}}$.
Consider the sets formed by the union of finitely many closed rays:
$${\cal C}_0 = {\cal C} \cup \overline{R_{t_1}} \cup \cdots 
\cup \overline{R_{t_k}},$$
$${\cal C}_1 = P^{-1} ({\cal C}_0);$$
and the corresponding simply connected domains containing the
Cremer fixed point $\hat{z}$: 
$$U_0 = U_{{\cal C}_0} ( \hat{z} ), $$
$$U_1 = U_{{\cal C}_1} ( \hat{z} ). $$

The rest of the proof consists in showing that $U_1$ contains at most
finitely many cycles. 
For each repelling cycle contained in $U_1$ there is at
least one cycle of external rays contained in $U_1$. Then, it is
enough to show that $U_1$ contains at most finitely many cycles
of external rays.\medskip

First observe that  the proper holomorphic map
$$P|_{U_1} : U_1 \rightarrow U_0$$
has no critical points.
Regular coverings of simply connected domains are
homeomorphisms. Therefore,  $P|_{U_1}$ is a conformal isomorphism.
\medskip

The polynomial $P$ acts in the angles of the external rays 
as the multiplication by $ d = \deg P$ map,
$$\begin{array}{lccc}
m_d : &  {\Bbb R}/{\Bbb Z} & \rightarrow & {\Bbb R}/{\Bbb Z}\\
      &      t             & \mapsto     &  d \cdot t \\ 
\end{array}$$

We are interested in the external rays contained in $U_1$ and
their image, the external rays contained in $U_0$. Denote
by $B_1$ and $B_0$ the corresponding sets formed by the angles
of such rays:
$$ B_0 = \{ t  :  R_t \subset {U}_{0} \} \subset {\Bbb R}/{\Bbb Z},$$
$$ B_1 = \{ t  :  R_t \subset {U}_{1} \} \subset {\Bbb R}/{\Bbb Z}.$$

As mentioned before the map $P|_{U_1} : U_1 \rightarrow U_0$
is a conformal isomorphism, therefore 
$$m_d|_{ B_1} : B_1 \rightarrow B_0 $$
is a cyclic order preserving homeomorphism.
Hence, $m_d|_{ B_1}$ admits a continuos extension
$g :  {\Bbb R}/{\Bbb Z}  \rightarrow  {\Bbb R}/{\Bbb Z}$
which is a  non-decreasing circle map (i.e. the lift is non-decreasing).
\medskip

A non-decreasing circle map $g$  either is free of periodic orbits or all
periodic orbits have the same period.  If $m_d$ has cycles contained
in $B_1$, then all of them have the same period.  The map 
$m_d : {\Bbb R}/{\Bbb Z} \rightarrow {\Bbb R}/{\Bbb Z}$ has finitely
many cycles of the same period.  We conclude that $m_d$ has at most
finitely many cycles contained in $B_1$. An equivalent statement is
that there are at most finitely many cycles of external rays contained
in $U_1$, which contradicts the assumption that the Cremer fixed point
$\hat{z}$ is approximated by small cycles. 
\medskip

Under that assumption that $\hat{z}$ is not a critical value,
we have proved that there exists a critical value $v \neq \hat{z}$ such that 
every ray that accumulates at $v$ also accumulates at $\hat{z}$.
It follows that $v$ is not accessible from ${\Bbb C} \smallsetminus K_P$
\medskip

In the case that $\hat{z}$ is a critical value
we outline how to modify the arguments used before.
Call $c_1, \dots , c_k$ the critical points of $P$ 
contained in $U_{\cal C} (\hat{z})$. List without 
repetition their critical values $v_0 = \hat{z}, v_1 ,\dots,
v_l$. The proof also proceeds by contradiction. This time, we assume  
that $R_{t_0}$ lands at $v_0$ and that $v_i \in \overline{R_{t_i}}$
and $\hat{z}  \notin \overline{R_{t_i}}$.
\medskip

 As before, consider the sets formed by the union of closed rays:
$${\cal C}_0 = {\cal C} \cup  \overline{R_{t_1}} \cup \cdots 
\cup \overline{R_{t_k}},$$
$${\cal C}_1 = P^{-1} ({\cal C}_0).$$
Now let $U_0$ (resp.$ U_1$) be the unique component of
the complement of ${\cal C}_0$ (resp. ${\cal C}_1$) such that
$\hat{z}$ lies in the boundary of $U_0$ (resp. $U_1$).
The same arguments used before show that $U_0$ contains at most finitely
many cycles. Since the union of $U_0$ and a preimage of 
$\overline{R_{t_0}}$ is 
an open neighborhood of $\hat{z}$, it follows that $\hat{z}$ cannot
be approximated by small cycles.
 \hfill $\Box$
\bigskip

\noindent
{\bf Theorem 1.1} {\it
Let $P$ be a polynomial with a Cremer fixed point $\hat{z}$
that is approximated by small
cycles. Then there exists a critical point which
is not accessible from ${\Bbb C} \smallsetminus J_P$.}
\medskip

For the case in which $J_P$ is connected, Theorem 1.1 is an immediate
consequence of Theorem 1.3.
In the case that $J_P$ is disconnected the idea is to extract from $P$
a {\it polynomial like map} which after {\it straightening} becomes a
Cremer polynomial with connected Julia set. Then we can apply Theorem
1.3  to obtain a non-accessible critical point. Before
going further we need some definitions and results about polynomial
like mappings (compare \cite{dh:p-l}):
\medskip

\noindent
{\bf Definition:} Let $D \subset {\Bbb C}$ and  $D^\prime \subset {\Bbb C}$
be bounded simply connected domains with smooth boundaries such that
$\overline{D} \subset D^\prime$. The triple $( f ; D , D^\prime )$ is 
called a {\it polynomial like map} of degree $d$, if the map
$f: \overline{D} \rightarrow \overline{D^\prime}$
is a $d$-fold branched covering ($d \geq 2$) which is holomorphic in $D$.
\medskip

The {\it filled Julia set} $K_{( f ; D , D^\prime )}$ of $( f ; D , D^\prime )$
is the set of points in $D$ for which the forward iterates of $f$ are
well defined:
$$K_{( f ; D , D^\prime )} = 
{\bigcap} f^{-n} ( \overline{D} ).$$
The {\it Julia set} $J_{( f ; D , D^\prime )}$ is the boundary of $K_{( f ; D , D^\prime )}$.

\medskip
Polynomial like mappings can be extended to the complex plane in such
a way that they are quasiconformally conjugate to a polynomial of the
same degree:

\begin{theo}[Douady and Hubbard]
If $ ( f ; D , D^\prime )$ is a polynomial like map of degree $d$,
then there exists a quasiconformal map $\phi : {\Bbb C} \rightarrow
{\Bbb C}$ and a polynomial $Q$ of degree $d$ such that $\phi \circ f =
Q \circ \phi$ on $\overline{D}$.  Moreover, $\phi (K_{( f ; D , D^\prime
)}) = K_Q$.
\end{theo}

Now we have all the ingredients to prove Theorem 1.1 except for the following:

\begin{lem}
Let $P$ be a polynomial and $C$ a connected component of $K_P$ such that
$P(C)=C$. If $C$ is not a repelling fixed point then there exists 
a polynomial like map $( P ; D , D^\prime )$ such that
 $C=K_{( P ; D , D^\prime )}$.
\end{lem}

Let us defer the proof of this Lemma.
\medskip

\noindent
{\bf Proof of Theorem 1.1}: 
Let $C$ be the connected component of $K_P$ that contains the 
Cremer fixed point $\hat{z}$. It follows that  $P(C)=C$.
Extract a polynomial like map  $( P ; D , D^\prime )$ as in Lemma 4.2.
After `straightening'   $( P ; D , D^\prime )$
we obtain a polynomial $Q$ and a homeomorphism $\phi$
such that  $\phi \circ f = Q \circ \phi$ on $\overline{D}$
and  $K_Q = \phi (C)$ (Theorem 4.1). Notice that $K_Q$ is connected.
Hence, $J_Q = \phi (\partial C)$ is connected. \medskip

Every neighborhood $U \subset D$ of $\hat{z}$ contains infinitely
many cycles of $P$, hence every neighborhood of 
 $\phi (\hat{z})$   contains infinitely many
cycles of $Q$. Therefore,  $\phi (\hat{z})$ is a Cremer fixed point of $Q$
that is approximated by small cycles.
Now we are under the assumptions of Theorem 1.3 for connected
Julia sets,  let $c \in J_Q$ be a critical
point of $Q$ that is not accessible from ${\Bbb C}- J_Q$.\medskip

It is not difficult to check that
 $\phi^{-1}(c)$ is critical point  of $P$ that belongs
to $\partial C \subset J_P$. Also,  paths in 
${\Bbb C} \smallsetminus J_Q$ correspond under $\phi^{-1}$ to
 paths in ${\Bbb C}\smallsetminus \partial C$.
Thus, the critical point $\phi^{-1} (c)$ of $P$ is not accessible
from  ${\Bbb C} \smallsetminus \partial C \supset{\Bbb C} \smallsetminus J_P$.
\hfill $\Box$
\bigskip

\newpage
In order to prove Lemma 4.2 we need to work with the 
{\it Green function} $G$ of $K_P$:

$$\begin{array}{lccl}
  G : &  {\Bbb C} & \rightarrow & {\Bbb R}_{\geq 0}\\
      & z & \mapsto &  {\lim}_{n \rightarrow \infty} 
\frac{\log_+ |P^{\circ n} (z)|}{d^n}   \\ 
\end{array}$$
where $\log_+ s = \max \{ 0 , \log s \}$.

It can be shown that $G$ is continuous, vanishes in $K_P$ and $G \circ P =
d \cdot G$. 
The Green function is harmonic in $ {\Bbb C}-K_P$, in particular
differentiable.
If $r > 0$ is a regular value of $G$, then  each component of 
$\{ z : G (z) \leq r \}$ is a closed topological disc with smooth boundary.
A point $z \in  {\Bbb C} \smallsetminus K_P$ is a critical point of $G$
if and only if 
some forward iterate of $z$ is a critical point of $P$.
Thus, $G$ has regular  values arbitrarily close to $0$.
(see Section 17 of \cite{milnor:lec})
\medskip

\noindent
{\bf Proof of Lemma 4.2:}
Let us denote by $V(r)$ the connected  component of $\{ z : G (z) \leq r \}$
containing $C$.
It is not difficult to check that:
$$ C = \bigcap V(r)$$
where the intersection is taken over
all the regular values $r$ of $G$.
\medskip

In particular, if a critical point $c$ of $P$ is such that $c \in
V(r)$ for all regular values $r$ of $G$, then $c \in C$. 
It follows that there exists a regular value $r_0$
of $G$ such that all the critical points of $P$ contained 
in $V(r_0)$ belong to $C$. 
\medskip

Consider the topological discs with smooth boundary 
$D^{\prime} = \mbox{int} V(r_0)$ and $D = \mbox{int} V(r_0 / d)$.
The proper holomorphic map $P: D \rightarrow D^{\prime}$ is a 
$d$-fold branched covering. Under the assumption that $C$ is not a 
repelling fixed point, $P|_D$ is not a conformal isomorphism.
Hence, $d \geq 2$ and  $(P; D, D^{\prime})$ is a polynomial
like map.\medskip

 The critical
points of the polynomial like map $(P; D, D^{\prime})$ 
are contained in $C$, therefore they 
never escape from $D$. It follows that the filled Julia set 
$K_{ (P; D, D^{\prime})}$ is connected
and  $C \subset K_{ (P; D, D^{\prime})}
\subset K_P$. 
 But $C$ is a connected component of $K_P$,
so in fact $C = K_{ (P; D, D^{\prime})}$. \hfill $\Box$
\bigskip

\noindent
{\bf Corollary 1.2} {\it The critical point of a quadratic polynomial $P$
with a Cremer periodic point $\hat{z}$ is not accessible from the complement of
the Julia set.}
\medskip

The proof of Corollary 1.2 consists in showing that Cremer periodic points
of quadratic polynomials are aproximated by small cycles.
This is obtained by putting together results of  Douady, Hubbard and
Yoccoz.
 According to Douady and  Hubbard a quadratic 
polynomial $P$ with a Cremer periodic point (not fixed)
is renormalizable. After an appropiate 
renormalization of $P$ we obtain a quadratic polynomial with a Cremer
fixed point which, according to  Yoccoz, is approximated by
small cycles:
\medskip

\noindent
{\bf Proof:} Let $n$ be the period of $\hat{z}$. It follows from the
work of Douady and  Hubbard \cite{dh:p-l} that there exists a
quadratic like map $( P^{\circ n} , D, D^{\prime})$ which 
after straightening becomes a quadratic polynomial $Q$ with
a Cremer fixed point. More precisely, call $\phi$ the quasiconformal
conjugacy between $( P^{\circ n} , D, D^{\prime})$ and $Q$ (Theorem 4.1),
 then $\phi (\hat{z})$ is a Cremer fixed point of $Q$. According
to  Yoccoz \cite{yo:linea},  $\phi (\hat{z})$ is approximated by small
cycles of $Q$. Therefore, $\hat{z}$ is approximated by small
cycles of $P^{\circ n}$. Now apply Theorem~1.1 to $P^{\circ n}$ and
obtain a non-accessible critical point of  $P^{\circ n}$.
Since the critical points of  $P^{\circ n}$ are preimages of
the critical point of $P$ and accessibility is a completely invariant
property, it follows that the critical point of $P$ is not accessible
from the complement of the Julia set.\hfill $\Box$

\end{document}